\documentstyle[12pt]{article}
\def\INF{\infty}

\def \Z{\hbox{$Z\hskip -5.2pt Z$}}

\def\sZ{\hbox{$\sc Z\hskip -4.2pt Z$}}

\def \R{\hbox{$I\hskip -3pt R$}}

\def \C{\hbox{$C\hskip -5pt \vrule height 6pt depth 0pt \hskip 6pt$}}

\def\qed{\hfill \hfill \ifhmode\unskip\nobreak\fi\ifmmode\ifinner
         \else\hskip5pt\fi\fi
 \hbox{\hskip5pt\vrule width4pt height6pt depth1.5pt\hskip 1 pt}}

\def\d{\delta}

\def\l{\lambda}

\def\o{\omega}

\def\Vir{\hbox{\bf\sl Vir}}

\def\sc{\scriptstyle}
\def\ssc{\scriptscriptstyle}
\def\dis{\displaystyle}
\def\OVER#1#2{\mbox{${{\ssc\,}\dis#1{\ssc\,}}
\over{{\ssc\,}\rb{-2pt}
{\mbox{$\dis#2$}}{\ssc\,}}$}}

\def\cl{\centerline}

\def\rar{\rightarrow}

\def\Rla{\Leftrightarrow}

\def\bs{\backslash}

\def\vs{\vspace*}
\def\vs#1{}
\def\rb{\raisebox}
\def\VS{\mbox{\rb{-7pt}{\,}}}
\def\ra{\rangle}
\def\la{\langle}
\def\ni{\noindent}
\def\hi{\hangindent}
\def\ha{\hangafter}

\def\cc{\hbox{${\ssc\,}c\hskip -4.5pt c{\ssc\,}$}}
\def\O{in Category $\cal O${}}

\begin{document}
\def\ABS
{It is proved that an indecomposable Harish-Chandra
module over the Virasoro
algebra must be (i) a uniformly bounded module, or (ii) a
module \O, or (iii) a module \O$^-$, or (iv) a module which contains
the trivial module as one of its composition factors.}
\def\KEYS
{Virasoro algebra, intermediate series, uniformly bounded, Category ${\cal O}$.}
\cl{\bf
ON INDECOMPOSABLE MODULES OVER THE VIRASORO ALGEBRA\footnote
{
This work is supported by a Fund from
National Education Department of China.}}
\par\
\par
\cl{Yucai Su}
\par\
\vs{-7pt}\par
\cl{Department of Applied Mathematics,
Shanghai Jiaotong University, Shanghai 200030, China
}
\par\
\par\
{\bf Abstract} \ABS
\par\ni
\par\
\par\
The Virasoro algebra \Vir, as the universal central extension of the
\vs{-2pt}infinite dimensional complex Lie algebra of the linear differential
operators $\{t^{i+1}{\OVER{d}{dt}}\,|\,i\in\Z\}$, of \vs{-2pt}interest to both
physicists and mathematicians [1-9], is the Lie algebra with basis\vs{-2pt}
$\{L_i,\cc\,|\,i\in\Z\}$
such that
$[L_i,L_j]=(j-i)L_{i+j}+{\OVER{i^3-i}{12}}\d_{i,-j}\cc,$
\vs{-2pt}$[L_i,\cc]=0,i,j\in\Z.$
\par
Set
$\Vir_\pm=\oplus_{i\in\sZ_\pm\bs\{0\}}\C L_i$, $\Vir_0=\C L_0\oplus\C \cc$,
$\Vir_{[i,j]}=\oplus_{i\le k\le j}\C L_k$, $i,j\in\Z$ and
$\Vir_{[i,\INF)}=\oplus_{i\le k}\C L_k$.
Then \Vir\, has the triangular decomposition
$\Vir=\Vir_-\oplus\Vir_0\oplus\Vir_+$
and we have the universal enveloping algebra decomposition
$U(\Vir)=U(\Vir_-)U(\Vir_0)U(\Vir_+)$.
\par
Consider a Harish-Chandra $\Vir$-module $V$, i.e., a module
with finite dimensional weight space
decomposition:
$$
V=\bigoplus_{\l\in\C}V_\l,\ V_\l=\{v\in V\,|\,L_0v=\l v\},
{\rm dim\sc\,}V_\l<\INF,\ \l\in\C.
\eqno(1)$$
Since the central element $\cc$ acts as a scalar on any indecomposable module,
we shall always suppose $\cc V=hV$ for some $h\in\C$ (and $h=0$ when $V$ is
uniformly bounded).
\par
Throughout this paper,
we shall always suppose that $V$ is a $\Vir$-module with decomposition (1).
\par
{\bf Definition 1}.
A module $V$ is said to be
(i) a {\it module of the intermediate series} [5] if it is indecomposable and
    ${\rm dim\sc\,}V_\l\le1$ for all $\l\in\C$;
(ii) a {\it uniformly bounded module} if there exists $N>0$ such that
     ${\rm dim\sc\,}V_\l\le N$ for all $\l\in\C$;
(iii) a {\it highest} or {\it lowest weight module} if it is generated by
      a highest or lowest weight vector $v$ with $L_iv=0$ for all $i>0$
      or $i<0$ respectively;
 (iv) a {\it module in Category $\cal O$} or {${\cal O}^-$} if there exists
      $\l_0\in\C$ such that for $\l\in\C$ with $V_\l\ne0$
      one has $\l\!\le\!\l_0$ or $\l_0\!\le\!\l$ respectively
      (here and below a partial order is defined on $\C$ by:
       $\l\!\le\!\mu\Rla\l\!-\!\mu\!\in\!\R_-$).
\par
A vector $v\ne0$ is said to be
    {\it primitive},
    {\it anti-primitive},
 {\it strongly primitive} or
 {\it strongly anti-primitive}
if $v\notin U(\Vir)\Vir_+v$, $v\notin U(\Vir)\Vir_-v$,
$\Vir_+v=0$ or $\Vir_-v=0$ respectively.
\qed\par
Highest weight modules are examples of modules \O. However
a module \O\, may have infinite number of composition factors
and may not be generated by any finite number of primitive vectors.
Simple modules over the Virasoro algebra have been studied in
[1,3,5-9]. In particular, we have
\par
{\bf Theorem 2}.
A simple module over the Virasoro algebra must be
  (i) a module of the intermediate series, or
 (ii) a highest weight module, or
(iii) a lowest weight module.
\qed\par
Thus simple modules are, in this sense, well known. This theorem was
conjectured in [3,4] and proved in [8] and also
partially in [1,7,9], and was later generalized to the
super-Virasoro
algebras in [10] and the higher rank Virasoro algebras (in some sense)
in [11].
\par
A module of the intermediate series [5] must be one of
$A_{a,b},A(a),B(a)$, $a,b\in\C$, or one of their
quotient submodules, where $A_{a,b},A(a),B(a)$ all
have a basis $\{x_k\,|\,k\in\Z\}$ such that $\cc$ acts trivially and
$$
\matrix{
\hfill A_{a,b}:&L_i x_k=(a+k+bi)x_{i+k},\VS\hfill&&\cr
A(a):&L_i x_k=(i+k)x_{i+k},\hfill&
\!\!\!\!\!\!\!\!\!\!\!\!\!\!\!\!\!k\ne0,\hfill&
L_i x_0&\!\!\!\!=
i(i+a)x_i,\VS\hfill\cr
B(a):&L_i x_k=k x_{i+k},\hfill&
\!\!\!\!\!\!\!\!\!\!\!\!\!\!\!\!\!k\ne-i,\hfill&
L_i x_{-i}&\!\!\!\!=-i(i+a)x_0,\hfill\cr
}
\eqno(2)$$
for $i,k\in\Z$. We have
$$
\mbox{(i) } A_{a,b}\mbox{ is simple }\Rla a\notin\Z
\mbox{ or }a\in\Z, b\ne0,1, \mbox{ and (ii) }A_{a,1}\cong
A_{a,0} \mbox{ if }a\notin\Z.
\eqno(3)$$
Our main result is the following.
\par
{\bf Theorem 3}.
An indecomposable module over the Virasoro algebra must be
  (i) a uniformly bounded module, or
 (ii) a module \O, or
(iii) a module \O$^-$, or
 (iv) a module which contains the trivial module $V(0)$ as one of its
      composition factors.
\par
Theorem 3 shows that if an indecomposable module does not contain
$V(0)$ as a composition factor, then its composition factors
are all of the same type:
  (i) modules of the intermediate series,
 (ii) highest weight modules,
(iii) lowest weight modules.
However, if an indecomposable module contains the composition factor $V(0)$,
then it can have all types of composition factors. As an example, let $V'$
be the Verma module $M(0)$ [2] generated by the highest weight vector $v'_0$
with weight $0$ such that $\cc v'_0{\sc\!}={\sc\!}0$,
let $V''$ be the anti-Verma module
$M_*(0)$ generated by the lowest weight vector $v''_0$ with weight $0$ such
that $\cc v''_0{\sc\!}={\sc\!}0$, let $V'''$ be
the module of the intermediate series of
type $A(a)$ in (2), and finally let $V$ be the submodule of
$V'{\sc\!}\oplus {\sc\!}V''{\sc\!}\oplus {\sc\!}V'''$ generated by
$v_0{\sc\!}={\sc\!}v'_0{\sc\!}+{\sc\!}v''_0{\sc\!}+{\sc\!}x_0$,
then $V$ is an indecomposable module which contains all types of composition
factors such that $V(0)$ is the {\it top} composition factor.
This example shows how one can get a module of type (iv) in Theorem 3
and how such a module can be like.
\par\
\par
\ni{\large\bf 1 \ Proof of Theorem 3}
\par
Let $V$ be an indecomposable module over \Vir. We can therefore suppose
$$
V=\sum_{k\in\sZ}V_{k+a}\mbox{ \ for some \ } a\in\C.
\eqno(4)$$
We always suppose that $0\le{\rm Re}(a)<1$.
\par
From now to Lemma 8, we suppose that $V$ is any module as in (4) such that
the trivial module $V(0)$ is not a composition factor of $V$.
\par
{\bf Lemma 4}.
Let $0\ne v\in V$, $i\in\Z_+$. If $\Vir_{[i,\INF)}v=0$, then the submodule
$V'=U(\Vir)v$ has
a highest weight vector (which is clearly strongly primitive in $V$).
\par
{\bf Proof.}
The proof is exactly similar to that of Lemma 3.1 in [10]
(see also Lemma 1.6 in [8]).
\qed\par
{\bf Lemma 5}.
For any $i\in\Z$, there exist only a finite number of primitive
(anti-primitive) vectors with weights $a+j$ such that $j\ge i$
(respectively $j\le i$).
\par
{\bf Proof.} Suppose there are infinite number of primitive vectors
$v_k$ with weights $\l_k=a+j_k$, $j_k>i$, $k=1,2,\cdots$. Let $U^{(k)}$ be
the composition factors corresponding to $v_k$.
Then
$${\rm dim\sc\,}V_{a+i-i_0}\ge\sum_{k\ge1}{\rm dim\sc\,}U^{(k)}_{a+i-i_0},
\mbox{ where } i_0=0\mbox{ or }1.
\eqno(5)$$
For any $k\ge1$, since $U^{(k)}\cong V(\l_k)$ is
a nontrivial simple highest weight $\Vir$-module
and $\l_k>a+i$, we always have either ${\rm dim\sc\,}U^{(k)}_{a+i}\ge1$
or ${\rm dim\sc\,}U^{(k)}_{a+i-1}\ge1$; if not, then
$L_{-1}U^{(k)}_{a+i+1}=L_{-2}U^{(k)}_{a+i+1}=0$ and since $\Vir_-$ is
generated by $L_{-1},L_{-2}$, we obtain that $\Vir_-U^{(k)}_{a+i+1}=0$
and that $U^{(k)}$ has a lowest weight,
this is impossible. This shows that the right-hand side of (5) is the
infinity for at least one $i_0=0$ or 1, this is
a contradiction to (1).
\qed\par
{\bf Lemma 6}.
Any submodule $U$ generated by a primitive vector
$v_\l$ of weight $\l$ is \O.
\par
{\bf Proof.}
Suppose conversely that $U$ is not \O. If $U$ has strongly primitive vectors,
let $W$ be sum of all submodules $U(v_\mu)$ of $U$ generated by
{\it strongly} primitive vectors $v_\mu$. Since each $U(v_\mu)$ is a highest
weight module, by Lemma 5, $W$ must be \O. Thus $U/W$ is still not \O.
For a given $k\in\Z_+$, since the number of the primitive vectors with
weights larger than $\l-k$ is finite, by considering the module $U/W$ instead
of $U$ and if necessary repeating the above process, we can suppose that
in $U/W$, there is no strongly primitive vector with weight $\mu>\l-k$, and
that $U/W$ is not \O. Now by Lemma 4, for any
$i\in\Z_+,\,0\ne x\in (U/W)_{\l-k}$, we have
$\Vir_{[i,\INF)}x\ne0$. Let $S$ be the Lie subalgebra of $\Vir$ generated
by $\Vir_{[k,k+1]}$, then clearly $S\supset\Vir_{[i,\INF)}$ for some $i$
(for any integer $m$ large enough, it can be written as $m=xk+y(k+1),
\,1\le x,y\in\Z$, so $L_m$ can be generated by $L_k, L_{k+1}$). Thus
$\Vir_{[k,k+1]}x\ne0$, i.e.,
${\rm ker}(L_k|_{(U/W)_{\l-k}})\cap{\rm ker}(L_{k+1}|_{(U/W)_{\l-k}})=0$.
That is,
$L_k|_{(U/W)_{\l-k}}\oplus L_{k+1}|_{(U/W)_{\l-k}}$: $(U/W)_{\l-k}\rar
(U/W)_\l\oplus (U/W)_{\l+1}$ is injective, and so ${\rm dim\sc\,}(U/W)_{\l-k}
\le {\rm dim\sc\,}(U/W)_{\l}+{\rm dim\sc\,}(U/W)_{\l+1}$. In particular,
$${\rm dim\sc\,}V(\l)_{\l-k}
\le {\rm dim\sc\,}(U/W)_{\l-k}
\le {\rm dim\sc\,}(U/W)_{\l}+{\rm dim\sc\,}(U/W)_{\l+1}
\le {\rm dim\sc\,}U_{\l}+{\rm dim\sc\,}U_{\l+1},
\eqno(6)$$
for $k\!\in\!\Z_+$, where $V(\l)$ is the simple highest weight module
with highest weight $\l$, which is the top composition factor of $U$.
The right-hand side of (6) is a fixed number.
This is impossible since every nontrivial simple highest weight
module $V(\l)$ is not uniformly bounded [1,3].
\qed\par
Exactly analogous to Lemma 6, we have
\par
{\bf Lemma 7}.
Any submodule $U'$ generated by an anti-primitive vector is \O$^-$.
\qed\par
{\bf Lemma 8}.
Let $U''$ be an indecomposable submodule generated by $x\in V$, where
$x$ corresponds to a basis element of $A_{a,b}$ in (2), such that
the top composition factor of $U''$ has type $A_{a,b}$, then $U''$
is uniformly bounded.
\par
{\bf Proof.} If $U''$ does not have primitive vectors or
anti-primitive vectors, then using similar arguments in the proof of
Lemma 6, we can show that $U''$ is uniformly bounded. Now suppose that
$U''$ has, say, a primitive vector.
Let $W',W''$ be the submodules of $U''$ generated by all primitive
vectors and by all anti-primitive vectors respectively. Then
by Lemmas 5-7,
$W'\cap W''=\{0\}$ since the only nonzero module both \O\ and
\O$^-$ is the trivial module.
Let $W=U''/W''$. Then $W$ is still indecomposable,
generated by $x$, having no anti-primitive vectors, and $W/W'$ has no
primitive or anti-primitive vectors, thus $W/W'$
is uniformly bounded. Hence by
[7], there exists $N\in\Z_+$ such that ${\rm dim\ssc\,}(W/W')_\l\le N$ for
all $\l\in a+\Z$.
\par
For any $\l\in a+\Z$, choose a basis $B_\l=C_\l\cup D_\l$ of $W_\l$ such that
$C_\l$ is a basis of $W'_\l$ and $D_\l$ corresponds to a basis of $(W/W')_\l$.
In the dual space $W^*_\l$ of $W_\l$, take the dual basis $B^*_\l=
C^*_\l\cup D^*_\l$ of $B_\l$. For any $v=\sum_{u\in B_\l}c_u u\in W_\l$,
$c_u\in\C$, writing as a linear
combination of elements of $B_\l$, we define the dual element $v^*=
\sum_{u^*\in B^*_\l} c_u u^*$, the same combination of elements
of $B^*_\l$.
Now in the dual space $W^*=\oplus_{\l\in a+\sZ}W^*_\l$ of $W$, we define a
$\Vir$-module structure as follows:
$$
\la L_i u^*_\mu,v_\nu\ra=\la u^*_\mu,L_{-i} v_\nu\ra,
\eqno(7)$$
for $i\in\Z,\,u_\mu\in W_\mu,\,v_\nu\in W_\nu,\,\mu,\nu\in a+\Z$.
\par
It is straightforward to verify that (7) defines a $\Vir$-module $W^*$.
\par
{\bf Claim 1}. For any $0\ne u^*_\l\in W^*_\l$, there exist
$i\in\Z_+\bs\{0\}$ and some $g^*\in U(\Vir)_i$ such that
$0\ne g^*u^*_\l\in W^*_{i+\l}$.
\par
Suppose $u^*_\l\ne0$, i.e.,
$u_\l=\sum_{u\in B_\l}c_u u$ with $c_{u_0}\ne0$ for some $u_0\in B_\l$.
Since $u_0$ is generated by $x$, which corresponds to a basis element of
the top composition factor of type $A_{a,b}$, there exist $\eta>\l$ and
some element $v_\eta\in W_\eta$ such that $u_0$ is also
generated by $v_\eta$,
i.e., there exists $g\in U(\Vir)_{-i}$, where $i=\eta-\l\in\Z_+\bs\{0\}$,
such that $u_0=g v_\eta$. Then
$$
0\ne c_{u_0}=\la u_\l^*,u_0\ra=
\la u^*_\l,g v_\eta\ra=\la \o(g)u^*_\l,v_\eta\ra,
\eqno(8)$$
where $\o$ is the anti-involution of $U(\Vir)$ defined by $\o(L_j)=L_{-j}$
for all $j\in\Z$, and where the last equality follows from definition (7).
This shows that $g^*u^*_\l\ne0$ for $g^*=\o(g)\in U(\Vir)_i$.
\par
Claim 1 in particular shows that $W^*$ has no strongly primitive vectors.
Then as in the proof of Lemma 6, $W^*$ is uniformly bounded.
This contradicts that $W$ is not uniformly bounded.
\qed\par
{\bf Proof of Theorem 3.}
Suppose that $V$ is not of type (iv). Let $W,W',W''$ be the submodules of $V$
generated respectively by all modules $U$, $U'$, $U''$ in Lemmas 6-8. Then
clearly, $W,W',W''$ are disjoint with each other and
$V=W\oplus W'\oplus W''$. Since $V$ is indecomposable, $V$ must be
one of $W,W',W''$.
\qed\par\
\par\
\cl{\large\bf References}
\par\ni\hi2ex\ha1
1 Chari V, Pressley A. Unitary representations of the Virasoro
 algebra and a conjecture of Kac. Compositio Math, 1988, 67: 315-342.
\par\ni\hi2ex\ha1
2 Feigin B L, Fuchs D B. Verma modules over the Virasoro algebra.
 Lecture Notes in Math, 1984, 1060: 230-245.
\par\ni\hi2ex\ha1
3 Kac V G. Some problems on infinite-dimensional Lie algebras and their
 representations. Lie algebras and related topics, Lecture Notes in Math,
 1982, 933: 117-126.
\par\ni\hi2ex\ha1
4 Kac V G. Infinite dimensional Lie algebras. 2nd ed,
 Birkh\"auser Boston, Cambridge, 1985.
\par\ni\hi2ex\ha1
5 Kaplansky I, Santharoubane L J. Harish-Chandra modules over the Virasoro
 algebra. Infinite-Dimensional Groups with Application,
 Math Sci Res Inst Publ, 1985, 4: 217-231.
\par\ni\hi2ex\ha1
6 Langlands R. On unitary representations of the Virasoro algebra.
 Infinite-Dimensional Lie algebras and Their Application,
 World Scientific, Singapore, 1986, 141-159.
\par\ni\hi2ex\ha1
7 Martin C, Piard A. Indecomposable modules over the Virasoro Lie algebra
 and a conjecture of V Kac. Comm Math Phys, 1991, 137: 109-132.
\par\ni\hi2ex\ha1
8 Mathieu O. Classification of Harish-Chandra modules over the
 Virasoro Lie algebra. Invent Math, 1992, 107: 225-234.
\par\ni\hi2ex\ha1
9 Su Y. A classification of indecomposable $sl_2(\C)$-modules
 and a conjecture of Kac on Irreducible modules over the Virasoro
 algebra. J Alg, 1993, 161: 33-46.
\par\ni\hi2ex\ha1
10 Su Y. Classification of Harish-Chandra modules over the super-Virasoro
 algebras. Comm Alg, 1995, 23: 3653-3675.
\par\ni\hi2ex\ha1
11 Su Y. Simple modules over the high rank Virasoro algebras. Comm Alg,
 in press.
\end{document}